\documentclass[aps,showpacs,preprint]{revtex4}
\usepackage{graphicx}
\usepackage{amssymb,amsmath}

\begin{document}

\title{The Riemann zeros and the cyclic Renormalization Group
}

\author{Germ\'an Sierra}

\affiliation{Instituto de F\'{\i}sica Te\'orica, CSIC-UAM, Madrid, Spain}
\date{October, 2005}

\date{October 2005}

\bigskip\bigskip\bigskip\bigskip

%
\font\numbers=cmss12
\font\upright=cmu10 scaled\magstep1
\def\stroke{\vrule height8pt width0.4pt depth-0.1pt}
\def\topfleck{\vrule height8pt width0.5pt depth-5.9pt}
\def\botfleck{\vrule height2pt width0.5pt depth0.1pt}
\def\Zmath{\vcenter{\hbox{\numbers\rlap{\rlap{Z}\kern
0.8pt\topfleck}\kern 2.2pt
                   \rlap Z\kern 6pt\botfleck\kern 1pt}}}
\def\Qmath{\vcenter{\hbox{\upright\rlap{\rlap{Q}\kern
                   3.8pt\stroke}\phantom{Q}}}}
\def\Nmath{\vcenter{\hbox{\upright\rlap{I}\kern 1.7pt N}}}
\def\Cmath{\vcenter{\hbox{\upright\rlap{\rlap{C}\kern
                   3.8pt\stroke}\phantom{C}}}}
\def\Rmath{\vcenter{\hbox{\upright\rlap{I}\kern 1.7pt R}}}
\def\Z{\ifmmode\Zmath\else$\Zmath$\fi}
\def\Q{\ifmmode\Qmath\else$\Qmath$\fi}
\def\N{\ifmmode\Nmath\else$\Nmath$\fi}
\def\C{\ifmmode\Cmath\else$\Cmath$\fi}
\def\R{\ifmmode\Rmath\else$\Rmath$\fi}

\begin{abstract}
We propose a consistent quantization of the Berry-Keating Hamiltonian
$ x p$,  which is currently discussed in connection with 
the non trivial zeros of the Riemann zeta function. 
The smooth part of the Riemann counting formula of the zeros
is reproduced exactly. The zeros appear, not
as eigenstates, but as missing states in the spectrum, 
in agreement with Connes adelic approach to the Riemann hypothesis. 
The model is exactly solvable and renormalizable, with
a cyclic Renormalization Group. These results 
are obtained by mapping the Berry-Keating model 
into the Russian doll model of superconductivity. Finally,  
we propose a generalization of these models in an attempt to 
explain  the oscillatory part of the Riemann's formula.  
\end{abstract}


\maketitle

\vskip 0.2cm

%
%
%
%
\def\oti{{\otimes}}
\def\lb{ \left[ }
\def\rb{ \right]  }
\def\tilde{\widetilde}
\def\bar{\overline}
\def\hat{\widehat}
\def\*{\star}
\def\[{\left[}
\def\]{\right]}
\def\({\left(}      \def\BL{\Bigr(}
\def\){\right)}     \def\BR{\Bigr)}
    \def\BBL{\lb}
    \def\BBR{\rb}
%
%
\def\zb{{\bar{z} }}
\def\zbar{{\bar{z} }}
\def\frac#1#2{{#1 \over #2}}
\def\inv#1{{1 \over #1}}
\def\half{{1 \over 2}}
\def\d{\partial}
\def\der#1{{\partial \over \partial #1}}
\def\dd#1#2{{\partial #1 \over \partial #2}}
\def\vev#1{\langle #1 \rangle}
\def\ket#1{ | #1 \rangle}
\def\rvac{\hbox{$\vert 0\rangle$}}
\def\lvac{\hbox{$\langle 0 \vert $}}
\def\2pi{\hbox{$2\pi i$}}
\def\e#1{{\rm e}^{^{\textstyle #1}}}
\def\grad#1{\,\nabla\!_{{#1}}\,}
\def\dsl{\raise.15ex\hbox{/}\kern-.57em\partial}
\def\Dsl{\,\raise.15ex\hbox{/}\mkern-.13.5mu D}
%
%
\def\ga{\gamma}     \def\Ga{\Gamma}
\def\be{\beta}
\def\al{\alpha}
\def\ep{\epsilon}
\def\vep{\varepsilon}
\def\dep{d}
\def\arc{{\rm Arctan}}
\def\la{\lambda}    \def\La{\Lambda}
\def\de{\delta}     \def\De{\Delta}
\def\om{\omega}     \def\Om{\Omega}
\def\sig{\sigma}    \def\Sig{\Sigma}
\def\vphi{\varphi}
\def\sign{{\rm sign}}
%
%
\def\CA{{\cal A}}   \def\CB{{\cal B}}   \def\CC{{\cal C}}
\def\CD{{\cal D}}   \def\CE{{\cal E}}   \def\CF{{\cal F}}
\def\CG{{\cal G}}   \def\CH{{\cal H}}   \def\CI{{\cal J}}
\def\CJ{{\cal J}}   \def\CK{{\cal K}}   \def\CL{{\cal L}}
\def\CM{{\cal M}}   \def\CN{{\cal N}}   \def\CO{{\cal O}}
\def\CP{{\cal P}}   \def\CQ{{\cal Q}}   \def\CR{{\cal R}}
\def\CS{{\cal S}}   \def\CT{{\cal T}}   \def\CU{{\cal U}}
\def\CV{{\cal V}}   \def\CW{{\cal W}}   \def\CX{{\cal X}}
\def\CY{{\cal Y}}   \def\CZ{{\cal Z}}

\def\rvac{\hbox{$\vert 0\rangle$}}
\def\lvac{\hbox{$\langle 0 \vert $}}
\def\comm#1#2{ \BBL\ #1\ ,\ #2 \BBR }
\def\2pi{\hbox{$2\pi i$}}
\def\e#1{{\rm e}^{^{\textstyle #1}}}
\def\grad#1{\,\nabla\!_{{#1}}\,}
\def\dsl{\raise.15ex\hbox{/}\kern-.57em\partial}
\def\Dsl{\,\raise.15ex\hbox{/}\mkern-.13.5mu D}
%
%
%
\font\numbers=cmss12
\font\upright=cmu10 scaled\magstep1
\def\stroke{\vrule height8pt width0.4pt depth-0.1pt}
\def\topfleck{\vrule height8pt width0.5pt depth-5.9pt}
\def\botfleck{\vrule height2pt width0.5pt depth0.1pt}
\def\Zmath{\vcenter{\hbox{\numbers\rlap{\rlap{Z}\kern
0.8pt\topfleck}\kern 2.2pt
                   \rlap Z\kern 6pt\botfleck\kern 1pt}}}
\def\Qmath{\vcenter{\hbox{\upright\rlap{\rlap{Q}\kern
                   3.8pt\stroke}\phantom{Q}}}}
\def\Nmath{\vcenter{\hbox{\upright\rlap{I}\kern 1.7pt N}}}
\def\Cmath{\vcenter{\hbox{\upright\rlap{\rlap{C}\kern
                   3.8pt\stroke}\phantom{C}}}}
\def\Rmath{\vcenter{\hbox{\upright\rlap{I}\kern 1.7pt R}}}
\def\Z{\ifmmode\Zmath\else$\Zmath$\fi}
\def\Q{\ifmmode\Qmath\else$\Qmath$\fi}
\def\N{\ifmmode\Nmath\else$\Nmath$\fi}
\def\C{\ifmmode\Cmath\else$\Cmath$\fi}
\def\R{\ifmmode\Rmath\else$\Rmath$\fi}

\def\barray{\begin{eqnarray}}
\def\earray{\end{eqnarray}}
\def\beq{\begin{equation}}
\def\eeq{\end{equation}}

\def\no{\noindent}

\def\gpar{g_\parallel}
\def\gperp{g_\perp}

\def\Jb{\bar{J}}
\def\dx{\frac{d^2 x}{2\pi}}

\def\rap{\beta}
\def\s{\sigma}
\def\spec{\zeta}
\def\comb{\frac{\rap\theta}{2\pi} }
\def\Ga{\Gamma}

\def\L{{\cal L}}
\def\g{{\bf g}}
\def\K{{\cal K}}
\def\I{{\cal I}}
\def\M{{\cal M}}
\def\F{{\cal F}}

\def\gpar{g_\parallel}
\def\gperp{g_\perp}
\def\Jb{\bar{J}}
\def\dx{\frac{d^2 x}{2\pi}}
\def\imag{\Im {\it m}}
\def\real{\Re {\it e}}
\def\Jbar{{\bar{J}}}
\def\kh{{\hat{k}}}
\def\Im{{\rm Im}}
\def\Re{{\rm Re}}

\section{Introduction}

The Riemann hypothesis (RH) is a central problem in Pure Mathematics
due to its connection with Number theory and 
other branches of Mathematics and Physics. 
The RH is the
statement that all the non trivial zeros
of the zeta function $\zeta(s)$ lie 
on the critical line $Re(s) = 1/2$ \cite{Edwards}. 
Most of the physical approaches to prove
the RH are inspired by the old Polya and Hilbert
conjecture, which states that the imaginary part, $E_a$,   
of the Riemann zeros $\zeta(\frac{1}{2} + i E_a)=0$, 
are the eigenvalues of a Hamiltonian and thus real 
numbers if the RH is true \cite{Watkins}. This approach is
supported by the statistical properties of the
zeros, the Montgomery-Odlyzko law \cite{Mont,Odl},  based on 
the Gaussian Unitary Ensemble distribution (GUE) 
\cite{Metha}, and also by the quantum chaos interpretation of the
oscillatory part of the Riemann counting formula of the zeros \cite{B-chaos}.
Berry has proposed that the 
Riemann dynamics is given by a classical chaotic 
Hamiltonian with isolated period orbits, whose 
quantization would yield a point like spectrum containing 
the $E_a's$ \cite{B-chaos}.
Other physical approaches are based on statistical 
mechanics 
\cite{Julia,BC}, superconformal invariance \cite{P}, 
supersymmetric quantum mechanics \cite{C}, etc 
(see \cite{Watkins,Elizalde,Rosu} for recent reviews). 
An interesting related approach is the construction
of a quantum mechanical potential containing in its 
spectrum the prime numbers \cite{Mussardo}. 

The starting point of our work is the Hamiltonian 
$H_{\rm BK} = x p$, proposed by Berry and Keating, which reproduces
at the semiclassical level 
the smooth part of the Riemann's formula giving the number
of non trivial zeros below a given number, $E$ \cite{BK1,BK2}.
This Hamiltonian is formal since 
a consistent quantization has not yet been found.
 In this paper we shall propose a solution of this problem 
defining $ x p$ on a lattice. 
This is  achieved working, not with   
 $x p$, but with its inverse $1/ x p$,
which turns out to be related to a class of QM models
with limit cycles, or rather centers, in the 
 Renormalization Group (RG). 

The idea that the RG may have limit
cycles was first considered by Wilson in 1971 
 \cite{Kwilson}, however
at the time no models with this behavior were known.
In the last few years limit cycles in the RG
have been discovered in various models in several physical
contexts, including nuclear physics \cite{nuclear},  
quantum field theory \cite{BLflow,LRS1}, quantum mechanics \cite{GW},
superconductivity \cite{RD1,RD2}, 
Bose-Einstein condensation \cite{Bose},
effective low energy  QCD \cite{QCD}, S-matrix models \cite{LRS2,LS}, 
few body systems and Efimov states \cite{nuclear,few-body}, etc
(for a review of see \cite{few-body}). 
The subject of duality cascades in supersymmetric gauge theory \cite{Klebanov}
is also suggestive of limit-cycle behavior. 
The possibility of chaotic flows has also been recently
considered \cite{GW,morozov}.

The relevant model for this work is 
the, so called, Russian doll BCS model
of superconductivity, and specially its QM version \cite{RD1,RD2}. 
In the RD model the standard pairing coupling $g$
flows periodically with the scale, a fact that is intimately
related to the existence of a series of bound states
of two electrons (Cooper pairs) 
whose energies scale as $e^{-n  \lambda} \; (n=0,1,2,\dots)$, 
where $\lambda$ is the period of the RG cycles. 
These bound states have a size that scales
as  $e^{n  \lambda} \; (n=0,1,2,\dots)$
which is the reason for calling them Russian dolls, 
by analogy with the popular matrioskas. 
The RG period is the main feature of a RG with
limit cycles. In the RD model it is given by
$\lambda = 2 \pi/h$, where $h$ is a coupling
in the Hamiltonian that breaks the time 
reversal symmetry.

Why should this model be related to the Hamiltonian
$x p$ which seems so far apart? Some hints lie on   
the following observations. In the RD model
the wave function of a Cooper pair with energy $E$ 
is given  approximately by $\psi(n) \sim 1/(n-E) ^{1 - ih}$,  
where $n=1,2, \dots, N$ label the energy levels
of the electron pairs. Since in this model
the energies $E_n$ converge towards zero, $E_n \sim e^{- n \lambda}$,  
there exist a bound state with $E = 0$, whose
wave function is  $\psi(n) \sim 1/n ^{1 - ih}$.
This form recalls the 
Dirichlet series of the zeta function,
$\zeta(s) = \sum_n 1/n^s$, with $s= 1 - ih$. 
The latter sum is the interaction term that is multiplied 
by the coupling $g$. Apparently, the RD model
should be connected to $\zeta(1 - i h)$ rather than to
$\zeta(1/2 - i E)$. 
It is worth to mention that 
in the cyclic sine Gordon model,
which also has a cyclic RG,  
the zeta function $\zeta(1-ih)$ appears
in the expression of the effective central
charge, where $h$ is related to
the period of the RG cycles as $\lambda = \pi/h$ 
\cite{LRS1,LRS2}.

The situation is different in the BK model
where the formal eigenfunctions of the normal
ordered Hamiltonian $(x p+ p x)/2$ are given 
by $\psi(x) = 1/x^{1/2 - i E}$, that resembles 
the Dirichlet series of  $\zeta(1/2 - i E)$. 
Compare to this, the RD wave function with $E=0$
has in the exponent of $n$  the factor 1 instead of 1/2, 
and the imaginary
part is fixed to a constant $h$ instead of being energy dependent. 
We shall show in this paper that these
two problems can be solved at once by relating the BK and
the RD models by means of a third model whose classical
Hamiltonian is simply the inverse of the BK one, namely
$1/x p$. In a deep sense, the BK and the RD models turn out to be the  
two faces of the same coin. 

The organization of the paper is as follows.
In section II we review the BK   
and the RD models. We also define 
the inverse model based on $1/xp$ 
and discuss its relation
with the previous ones.  
In section III we solve
the inverse model in the continuum limit. In section IV
we study its renormalization showing the existence
of RG cycles and their connection with the
continuous solution. In section V we give its 
exact solution,  derive
the smooth part of the Riemann counting formula
and present a numerical analysis. 
In section VI we generalize the  inverse model
and begin to study its properties. We state the conclusions
and prospects in section VII.

\section{The Hamiltonians}

In this section we introduce the Hamiltonians
of the Berry-Keating model, the inverse model ($\I$)
and the Russian doll  model,  and study their relations
which can be represented symbolically as

\beq
BK \rightarrow \I = BK^{-1} \rightarrow RD 
\label{symbol}
\eeq

\no 
Before presenting the details, 
we shall start with an overview.

\no 
As we explained above, the  classical version
of the   Berry-Keating Hamiltonian $H_{BK}$
is  the product
$ x p$, where $x$ and $p$ are the position and momentum
of a particle moving in one dimension. $H_{BK}$
fullfills some of the asumptions of the quantum chaos approach
to the RH, particularly the breaking of time-reversal symmetry, 
to accomodate the 
GUE hypothesis and suggestive analogies involving the
zeta function, the trace formula, the Riemann-Siegel formula, etc. 
The BK proposal 
has however remained at a speculative level since a consistent QM
 model of $xp $ has not been constructed in so far. 

To solve this  problem we shall define 
the Hamiltonian $H_{I}$, which at classical
level is the inverse of the BK Hamiltonian, namely
$1/x p$. It would seem that nothing is gained by
this trick. Nevertheless, a consistent
QM model can be constructed by quantizing $1/x p$
on a lattice. 
The key observation is that the inverse
of the momentum operator $p = - i \hbar \; d/ dx$
is the 1D Green function 
$\frac{i}{\hbar} G(x, x')=  
\frac{i}{2 \hbar} \sign(x-x')$, where 
 $\sign(x)$ is the sign function. 
We shall define a regularized lattice version
of the operator $p^{-1}$ as  the matrix
$P^{-1}_{n,m}= \frac{i}{2 \hbar} \;  \sign(n-m)$, where $n$ and $m$
run over the integers $1,2, \dots, N$. The quantization
of the position operator $x$ is 
the diagonal matrix $X_{n,m} = n \delta_{n,m}$. In this construction
the Hamiltonian $H_I$ acts on a discrete Hilbert space 
of dimension $N$. 

It will be  important to study the behaviour of 
 $H_I$ under the Renormalization Group
transformation which integrates the highest energy level
reducing the size of the system to $N-1$.  
For the model to be renormalizable one needs an extra
term with  coupling constant $g_I$. The other coupling constant,
$h_I = 1/\hbar$, multiplies 
the term containing $P^{-1}_{n,m}$. Hence the $\I$ model
depends on two couplings constants $h_I$ and $g_I$. 
Under the RG transformations
$h_I$ remains invariant, while $g_I$ flows periodically with
the scale. 
This periodicity is in turn related to 
the spectrum of the $\I$ model.

The third Hamiltonian is the QM
Russian doll Hamiltonian, 
which is a generalization of the
standard BCS model of superconductivity having
a cyclic RG. The RD Hamiltonian $H_{RD}$ 
also depends on two coupling constants
$g_D$ and $h_D$, where $h_D$ multiplies 
a time-reversal breaking term, proportional 
to   $P^{-1}_{n,m}$, while
$g_D$ multiplies the familiar $s$-wave pairing interaction 
of the BCS model which preserves the time-reversal
symmetry. The spectrum of this model
has a series of bound and antibound states
with Russian doll scaling.
The RD model is brought in  due to its close
relationship with the $\I$ model. Indeed, we shall
show  that each eigenstate of the $\I$ model
coincides with the zero energy bound state of an
associated RD model. This bound state appears
at the threshold of that model. From the RG viewpoint
these two models will also be related.

\subsection{The Berry-Keating model}

The classical
BK Hamiltonian  \cite{BK1,BK2}

\beq
H^{\rm cl}_{BK} = x p,  
\label{a1}
\eeq

\no 
has classical trayectories 
given by the  hyperbolas

\beq
x(t) = x_0 \; e^{t} , \quad p(t) = p_0 \;  e^{-t}.  
\label{a2}
\eeq

\no
The dynamics is unbounded and 
one should not expect a discrete spectrum at the
quantum level. Despite of that, Berry and Keating 
regularized the model introducing the Planck cell
in phase space with sides $l_x$,  $l_p$ and area
$h = l_x l_p$, such that $|x| > l_x$ and $|p| > l_p$ \cite{BK1}. 
They computed the number
of states $N_{BK}(E)$, with an energy below $E$, using the
semiclassical formula

\beq
N_{BK}(E) = \frac{A(E)}{h} = \frac{1}{h} 
\int_{l_x}^{E/l_p} d x \; \int_{l_p}^{E/x} dp + \dots 
\label{a3}
\eeq

\no
where $A(E)$ is the area of the 
region $l_x < x < E/l_p, \;  l_p < p < E/l_x, \;  x p < E$ (see fig.1a). 
The result is

\beq
N_{BK}(E) = \frac{E}{h} \left( \log \frac{E}{h} -1 \right) + 1 + \dots
\label{a4} 
\eeq

\no
Taking into account a Maslov phase $-1/8$, BK finally obtain

\beq
N_{BK}(E) = \frac{E}{2 \pi} \left( \log \frac{E}{2 \pi} -1 \right) + 
\frac{7}{8} + \dots 
\label{a5} 
\eeq

\no
where the energy $E$ is measured in units of $\hbar$. 
Eq.(\ref{a5}) agrees with the asymptotic 
expansion of  the smooth part of the Riemann 
counting formula of the non trivial zeros of the zeta function

\beq
N_{sm}(E) = \frac{1}{\pi}  \Im \log \Gamma \left( \frac{1}{4} + 
\frac{i}{2} E \right) 
- \frac{E}{2 \pi} \log \pi + 1. 
\label{a6} 
\eeq

\no 
This expression is $\pi^{-1}$ times the phase of the
zeta function $\zeta(1/2 + i E)$ \cite{Bha}. 
The non trivial zeros, $\rho$, are known
to lie in the critical strip 
$ 0 < Re(\rho) <  1$. The number of them, $N(E)$,  in the range 
$ 0 < \Im(\rho) < E, $  is given by the Riemann formula, 

\beq
N(E) = N_{sm}(E) + N_{\rm osc}(E),  
\label{a7} 
\eeq

\no
where the oscillatory part depends of the zeta function
on the critical line,

\beq
N_{\rm osc}(E) = \frac{1}{\pi} \Im \log \zeta \left( 
\frac{1}{2} + i E \right) 
\label{a8} 
\eeq

\no
and it is of order $\log E$ \cite{Connes}.  
The derivation of eq.(\ref{a5}) is heuristic 
and, to our knowledge, has not been 
reproduced from the
quantization of the BK Hamiltonian.
 Bhaduri et al. \cite{Bha}  
obtained eq.(\ref{a6}), up to some factor, studying the 
inverted harmonic oscillator
Hamiltonian $p^2 - x^2$, which is related
by a canonical transformation
to $x p$.  They employed a phase
shift approach instead of the semiclassical 
bound state approach.

The BK work was motivated by Connes's work  
to prove the RH using p-adic numbers.
This is known as the adelic approach to the RH
and it is based on the construction of an abstract space  
where acts an hermitean operator, whose eigenvalues are
the non trivial zeros of the zeta function \cite{Connes}. 
The truth of the RH lies
in the proof of a certain classical trace formula.  
Connes has also considered the operator $x p$
in the adelic theory. The regularization is not
the Planckian cell, but a
standard cutoff in phase space given by 
$|x| < \Lambda$, $|p| < \Lambda$. The number of states 
in the region $0 < x < \Lambda, 0 < p < \Lambda, x p < E$
(see fig.1b)

\beq
N_{Co}(E) = \frac{1}{h}
\left[ 2 E -  \left(  \frac{E}{\Lambda} \right)^2 +    
\int_{E/\Lambda}^{\Lambda} d x \; \int_{E/\Lambda}^{E/x} dp + \dots 
\right]
\label{a9}
\eeq

\no
yields,

\beq
N_{Co}(E) = \frac{E}{2 \pi} \log \Lambda^2
- \frac{E}{2 \pi} \left( \log \frac{E}{2 \pi} -1 \right) + \dots
\label{a10} 
\eeq

\begin{figure}[t!]
\begin{center}
\includegraphics[height= 6 cm,angle= 0]{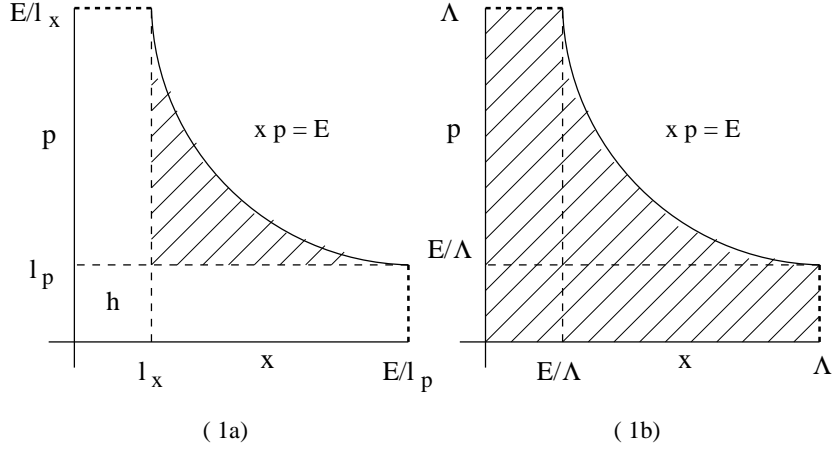}
\end{center}
\caption{
1a) The region in shadow is the phase space considered
by Berry and Keating in the semiclassical computation of the number
of states of the $xp$ Hamiltonian (see eq.(\ref{a3})).
1b) The shadow region gives the phase space considered
by Connes (see eq.(\ref{a9})).  
}
\label{berry}
\end{figure}

\no
where $E$ is measured in units of $\hbar$. 
There are important differences between eqs.(\ref{a10})
and (\ref{a5}). In Connes's formula
the number of states diverges with the cutoff. 
There is also a  negative sign of the
term in (\ref{a10}), associated to the asymptotic 
Riemann formula, 
as compared with the  BK formula (\ref{a5}). 
According to Connes, 
this negative sign agrees with the overall 
sign in a formal expression of 
$N_{\rm osc}(E)$ obtained using the Euler product
formula of $\zeta(1/2 + i E)$ into (\ref{a8}) \cite{Connes}. 
These results imply that   
the Riemann zeros, in the adelic approach,   
are missing states
in a continuum spectrum. The physical picture is
that of white light with dark absortion lines
labelled by the zeros of the zeta function. 
Using this analogy, 
the Riemann zeros, according to BK, would be emission lines. 
Our results are along the lines of Connes, as we shall show.

We shall end this brief review of the BK approach by
giving the formal hermitian operator
which corresponds to (\ref{a1}),

\beq
H_{BK} = \frac{1}{2} (x p + p x) =
- i \hbar \left( x  \frac{d}{dx} + \frac{1}{2} \right).    
\label{a11}
\eeq

\no
whose formal eigenfunctions are

\beq
\psi_E(x) = \frac{A}{x^{ 1/2 - i E/\hbar}}.  
\label{a13}
\eeq

\no
As noticed by BK, the wave function contains the power
$x^{-s}$ which appears in the Dirichlet series for $\zeta(s)$
and the Euler product, although those formulas
are convergent only in the region $Re(s) > 1$, while
in eq.(\ref{a13}) $s = 1/2 - i E/\hbar$,  which lies outside.  
The rest of the BK paper discusses
possible quantum boundary conditions 
that could generate  the Riemann zeros. The guiding idea
is that $x p$ is the generator of the scale transformations
in phase space, a fact that also plays a role in Connes's work.

{\bf Generalization of the BK model}

The  Hamiltonian (\ref{a1}) can be generalized 
to the form 
$ v(x) p$, where  $v(x)$ is a generic function which
we shall assume is positive and monotically increasing 
for $x >0$. 
The classical evolution equation $\dot{x}= v(x)$, implies that 
$v(x)$ is the velocity of the particle. The semiclassical
formula for the number of states \`a la Connes gives

\beq
N_{Co}(E) =  \frac{E}{2 \pi} \int_{v^{-1}(E/\Lambda)}^\Lambda
\frac{ dx }{v(x)} + \frac{\Lambda}{2 \pi} v^{-1}(E/\Lambda), 
\label{a14}
\eeq

\no
where $v^{-1}$ is the inverse function of $v(x)$.
If $v(x) = x$,  eq.(\ref{a14}) becomes
eq.(\ref{a10}). Similarly,  the formal Hamiltonian
associated to $v p$

\beq
H_{BK} = \frac{1}{2} (v(x) \;  p + p \;  v(x)) =
- i \hbar \left( v(x)  \frac{d}{dx} + \frac{1}{2} v'(x) \right),    
\label{a15}
\eeq

\no
has formal eigenfunctions 

\beq
\psi_E(x) = \frac{A}{v(x)^{1/2}}
{\rm exp} \left( \frac{ i E}{\hbar}
\int_{x_0}^x \; \frac{ d x'}{v(x')} \right).  
\label{a16}
\eeq

\no
If  $v(x) = x$ we recover  (\ref{a13}).

\subsection{The inverse Hamiltonian}

The inverse of the generalized BK Hamiltonian (\ref{a15}) is

\beq
H_{BK}^{-1} = v^{-1/2} \; p^{-1} \; v^{-1/2}, 
\label{a17}
\eeq

\no 
where we use the normal ordering 
$v^{1/2} \; p \; v^{1/2}$ of $H_{BK}$, which
is equivalent to (\ref{a15}).
By asumption $v(x) > 0$ for $x > 0$, then
the square root of $v(x)$ is real. For this
and other reasons 
will shall work on the half line $ x > 0$. 
The inverse of the momentum operator $p$ is

\beq
\langle x| p^{-1} |x' \rangle = 
\frac{i}{\hbar} G(x,x') = \frac{i}{2 \hbar} \; 
\sign(x-x'),  
\label{a18}
\eeq

\no
where $G(x,x')$ is the 1D Green function associated to
$d/dx$ and $\sign(x)$ is the sign function. 
The advantage of using $H_{BK}^{-1}$ instead of $H_{BK}$
is that the former operator can be easily regularized
on a 1D lattice due to the simplicity of 
the 1D Green function. The lattice
version of the formal Hamiltonian (\ref{a17}) 
is  the $N$-dimensional matrix

\beq
\langle n| H_I | m \rangle
= \frac{1}{2}  f_n \; \left( g_I + i h_I \sign(n-m) \right) f_m, 
\label{a19}
\eeq

\no
where $n,m = 1, 2 \dots, N$ and

\beq
h_I = \frac{1}{\hbar}, \qquad f_n = \frac{1}{\sqrt{v_n}}.  
\label{a20}
\eeq

\no
$H_I$ contains an extra term with coupling constant
$g_I$, which is needed for renormalization 
(see section IV). 
The {\em pure} BK model corresponds to the case $g_I= 0$, 
but we shall also consider non vanishing values of $g_I$.  
The regularized BK Hamiltonian will then be defined
as the inverse of $H_I$, namely

\beq
H_{BK}^{{\rm (reg)}}(v_n,g_I,\hbar) = H_I^{-1}(f_n,g_I,h_I)  
\label{a21-b}
\eeq

\no
In the case where $g_I =0$, one gets

\beq
H_{BK}^{{\rm (reg)}}(v_n, g_I = 0, \hbar) 
= v(X)^{1/2} \; P \;  v(X)^{1/2},  
\label{a21}
\eeq

\no
where $X$  and $P^{-1}$   are the $N$ dimensional matrices

\beq
\langle n| X | m \rangle = n \delta_{n,m}, \qquad 
\langle n| P^{-1} | m \rangle = \frac{i}{2 \hbar} \sign(n-m).  
\label{a22}
\eeq

\no
For $N$ an even number, the matrix $P^{-1}$ is not singular
and its inverse is given by

\beq
\langle n| P | m \rangle = - 2 \hbar i  \;  (-1)^{n+m} \;  \sign(n-m).  
\label{a22-2}
\eeq

\begin{figure}[t!]
\begin{center}
\includegraphics[height= 6 cm,angle= 0]{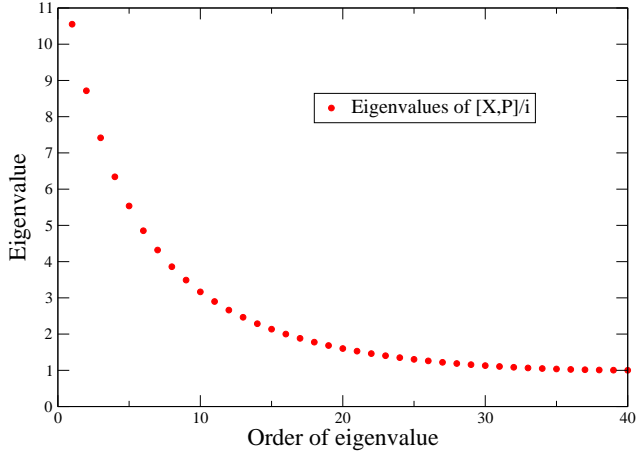}
\end{center}
\caption{
Eigenvalues of the commutator $[X,P]/i \;(\hbar =1)$ 
for $N=50$. We display 40 eigenvalues which converge towards
1. The remaining ones are larger in absolute value. 
The sum of all the eigenvalues is zero as 
a consequence of the vanishing
of the trace of $[X,P]$. 
}
\label{commuta}
\end{figure}

\no
$X$ and $P$ do not satisfy 
the canonical commutation relation $[x,p] = i \hbar$. 
The trace of the LHS of the commutator 
is zero while on the RHS
is $i \hbar N$. This is the textbook argument
to show that the Heisenberg's  indetermination relation
cannot be realized by finite dimensional matrices. 
Nevertheless, most of the  eigenvalues of
 $[X,P]/i$ converge towards $\hbar$ (see fig.\ref{commuta}). 
From this result we expect to recover the formal BK 
model in the continuum limit of the discrete 
Hamiltonian (\ref{a21}).
Since  $[X,P] \neq i \hbar$, the Hamiltonian
(\ref{a21}) cannot be written 
as $1/2( v(X) P + P v(X))$, but this fact is unimportant. 
We can then work either with the 
Hamiltonians (\ref{a19}) or (\ref{a21}).
We shall choose  the former for convenience.

The $g_I$ interaction has a classical version. This 
term  is non local in $x$-space and 
ultralocal in $p$-space. By dimensional reasons
it should be proportional to  $\delta(p)$,  yielding
$H_I^{\rm cl} = (1/ p + \delta(p)  g_I/2)/v(x)$, whose
inverse is again $p v(x)$ (recall that $p \; \delta(p) =0)$. 
At the quantum level $g_I$ does play a role
in  $H_{BK}^{\rm}(g_I) = H_I^{-1}(g_I)$ setting the boundary
conditions.

The operator $K= h/xp$
was used by BK to implement a canonical transformation
$ x \rightarrow x_1 = h/p, \; p \rightarrow p_1 = x p^2/h$, 
named {\em quantum exchange} by the $h$ dependence \cite{BK1}. 
Berry and Keating, tried to combine the dilatation
and the quantum exchange symmetries, 
to generate the Riemann zeros.
In our case, the operators $x p$ and $1/x p$ are used
as Hamiltonians and not as symmetries.

$H_I$ has two properties which
simplifies its study: 
it is renormalizable and exactly solvable for all
values of the coupling constants $h_I, g_I$
 and $f_n$. One can 
find an analytic form of the exact eigenvectors and
an explicit equation for 
the eigenenergies. These properties
are reminiscent of the Bethe ansatz, as it is indeed
the case.

\subsection{The Russian doll Hamiltonian}

A Hamiltonian closely related to $H_I$ is \cite{RD1}

\beq
\langle n| H_{RD} | m \rangle
= \vep_{n} \;  \delta_{n,m}
- \frac{1}{2}  \; \left( g_D + i h_D \;  \sign(n-m) \right),  
\label{a23}
\eeq

\no
where $\vep_n, g_D, h_D$ are  real parameters.
In physical applications $\vep_n$,
are the kinetic energies of pairs of electrons 
occupying doubly degenerate levels 
in the conduction band of a metal. 
The term proportional to $g_D$ is a pairing interaction
due to the phonon exchange in the $s$-wave channel,
which leads to the formation of Cooper pairs if $g_D  > 0$. 
The model with $h_D=0$ was first 
considered by Cooper and it was the  precursor
of the Bardeen-Cooper-Schrieffer model of superconductivity 
\cite{BCS,BCS-book}. 

\begin{figure}[t!]
\begin{center}
\includegraphics[height= 5 cm,angle= 0]{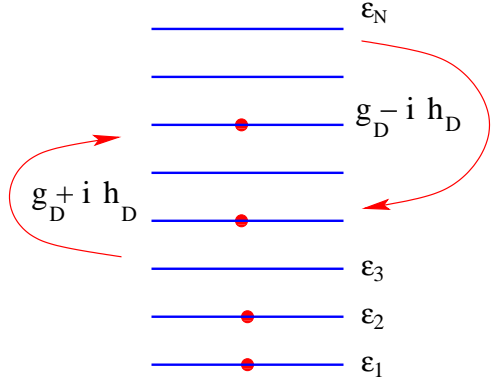}
\end{center}
\caption{
Pictorial representation of the RD model. The horizontal
lines represent the energy levels $\vep_n \;(n=1, \dots, N)$.
The dots are pairs of electrons occupying a doubly degenerate
state. The arrows give the transitions among levels
induced by the interactions in the Hamiltonian (\ref{a23}).  
}
\label{levelhopping2}
\end{figure}

The coupling $h_D$ was introduced
in reference \cite{RD1}, motivated by the 
work of Glazek and Wilson,  who proposed an 
extension of the well known 2D delta
function potential in order to 
show the existence of 
RG limit cycles in a simple QM model
\cite{GW}. 
The QM Russian doll model (\ref{a23}), 
together with its many body version, also have a cyclic RG, 
which has been studied in detail in references \cite{RD1,RD2}.
The exact solution was obtained in \cite{links}
using algebraic Bethe ansatz methods, showing
that the RD model is nothing but an inhomogenous 
 XXX vertex model with a boundary operator. 
In the vertex formulation the coupling $h_D$
parameterizes the quantum Yang-Baxter matrix, 
the coupling $g_D$ parameterizes the boundary operator and 
the energy levels $\vep_n$ give the inhomogenities. 
The most natural choice of the energy
levels is $\vep_n = n$,  in units of twice the  
electronic level spacing.
The equally spaced BCS model
for a finite and small number of energy levels
have been studied intensively due to the
fabrication of ultrasmall metallic grains
(see \cite{vDR} for a review). This model
is also relevant in Nuclear
Physics, where it is known as the picket fence model,
and in other potential applications in
Quantum Optics, and dilute Fermi-Bose gases
 (see \cite{duke} for a review). 
Up to date it is not known wether the RD-BCS model
of superconductivity explains some physical system 
existing in Nature or 
in the Lab. A possible reason 
is that the $h_D$ coupling breaks the time reversal
symmetry, while the microscopic laws of Nature
do not. However, this symmetry can
be broken spontaneously  or explicitely
by the action, for example,  of magnetic fields.

We shall next show the relation between the RD 
and the $\I$ models. Let's first compare their  
Schr\"odinger equations,

\beq
H_I \;  \psi = E_I \;  \psi, \qquad H_{RD} \;  \phi = E_{RD} \; \phi,
\label{a24}
\eeq

\no
which read explicitely

\barray
& E_I \psi_n =  \frac{1}{2} 
\sum_{m=1}^N f_n ( g_I + i h_I \;  \sign(n-m) ) f_m \psi_m  & 
\label{a25} \\
& (\vep_n - E_{RD})  \phi_n =  \frac{1}{2} 
\sum_{m=1}^N  ( g_D + i h_D \;  \sign(n-m) ) \phi_m. &  
\label{a26} 
\earray

\no
Given an eigenstate, $\psi$, of the $\I$ model let us define
the wave function

\beq
\phi_n = f_n \psi_n,  
\label{a27} 
\eeq

\no
writting (\ref{a25}) as

\beq
\frac{1}{f_n^2}  E_{I}  \phi_n =  \frac{1}{2} 
\sum_{m=1}^N  ( g_I + i h_I \;  \sign(n-m) ) \phi_m.   
\label{a28} 
\eeq

\no If $E_I \neq 0$, and dividing  eq.(\ref{a28})
by  $E_I$, one obtains eq.(\ref{a26})
with $E_{RD} =0$. This leads to the 
identifications:

\beq
\vep_n = \frac{1}{f_n^2} , \quad g_D = \frac{g_I}{E_I}, 
\quad h_D = \frac{h_I}{E_I},  
\label{a29} 
\eeq

\no which  establish a map from the spectrum
of the $\I$ model into a collection of RD models
which have in common the energy values, $\vep_n$,  
focusing on the $E_{RD}=0$ state, while varying
the couplings $g_D$ and $h_D$ according to $E_I$. 
Moreover, eqs.(\ref{a20}) establish a link 
between the regularized BK model and the RD model,
Combining (\ref{a29}) with (\ref{a20}) we get

\beq
\vep_n = v_n,  
\quad h_D = \frac{E}{\hbar},  
\label{a30} 
\eeq

\no where $E = 1/E_I$ is an eigenvalue
of $H_{\rm BK}= H_I^{-1}$. Let us remark that 
the transformation (\ref{a27})
is non unitary. For the model $\vep_n = n$ and in the 
limit $N \rightarrow \infty$,  
$\phi_n$ will be a normalizable wave function,
corresponding to a bound state,  
but $\psi_n$ will not, as corresponds to a scattering state. 

The map BK $ \rightarrow$ RD is quite remarkable. 
It  implies that the spectrum of the 
BK Hamiltonian can be found by looking at the zero
energy bound states of a RD model where the $h_D$ coupling
is fine tuned according to the energy of the 
state.  In this mapping, the velocities $v_n$ 
become energy levels $\vep_n$. In particular, the  choice
$v=x$ corresponds to the equally space model
$\vep_n = n$, which is the common choice 
for physical applications \cite{vDR,duke}. 
 
The correspondence (\ref{a30})
can be easily derived at the classical
level. Take  $g_D =0$ for simplicity. The classical energy 
(\ref{a23}) is $E_{RD} = \vep(x) - h_D/ p$
(with $\hbar = 1$). 
Setting $E_{RD} =0$ one finds $p \;  \vep(x) = h_D$,
which is the classical energy 
of the BK model with $v(x) = \vep(x)$ and 
$E = h_D$. 

In summary,  we have constructed a consistent quantum
model which generalizes the BK Hamiltonian $xp$
by means of the inverse  model $I$, whose
spectrum can be mapped into the zero eigenstates
of associated RD models. As we shall see this connection is the
key of the renormalizability and 
solvability of the BK model so defined.

\section{Continuum limit}

The three models described in the previous section can
be solved exactly. However it is worth to 
solve them first in the limit where the lattice size 
$N$ is very large. In this limit the discrete
variable $n$ will be considered as continuous 
and varying in the interval
$ 1 \leq n \leq N$. The Schr\"odinger eq.(\ref{a25})
becomes,

\beq
E_I \psi(n) =  \frac{1}{2} 
\int_{1}^N dm \;  f(n) ( g_I + i h_I \;  \sign(n-m) ) f(m) \psi(m),  
\label{b1}
\eeq

\no
where $\psi(n)$ and $f(n)$ are continuous functions. 
To solve (\ref{b1})  we use again the change of variables
(\ref{a27}),

\beq
\vep(n)  \phi(n) =  \frac{1}{2} 
\int_{1}^N  dm \;  ( g_D + i h_D \;  \sign(n-m) ) \phi(m),  
\label{b2} 
\eeq

\no
where $\vep(n), g_D$ and $h_D$ are given by eqs.(\ref{a29}).
As shown  above, eq.(\ref{b2}) is satisfied by 
an eigenstate of the RD model with $E_{RD}=0$, and can be solved
in the same way as was done in references \cite{RD1,RD2}.  
Taking the derivative with respect to $n$ 
yields

\beq
\frac{d}{d n} \left[ \vep(n)  \phi(n) \right]  =  i h_D \;  \phi(n),   
\label{b3} 
\eeq

\no
whose integral determines the functional
form of $\phi(n)$,

\beq
\phi(n) = \frac{A}{\vep(n)} 
{\rm exp} \left( i h_D \int_1^n \frac{d n'}{\vep(n')} \right).  
\label{b4} 
\eeq

\no
Using eqs.(\ref{a29}) and (\ref{a30}) this leads to

\beq
\psi_E(n) = \frac{A}{\vep(n)^{1/2}} 
{\rm exp} \left( \frac{i E}{ \hbar}  \int_1^n \frac{d n'}{\vep(n')} \right).  
\label{b5} 
\eeq

\no
which coincides with the wave function (\ref{a16}) 
for the generalized BK model. This result
agrees with the fact 
that  the commutator $[X,P]$ converges
asymptotically to $i \hbar$ (see fig.\ref{commuta}).  
To find the eigenenergies $E= 1/E_I$
let us consider eq.(\ref{b2}) at the boundaries
of the interval, $n=1, N$

\barray 
\vep(1)  \phi(1) & = &  \frac{1}{2}  ( g_D - i h_D)  
\int_{1}^N  dn \;   \phi(n)  
\label{b6}  \\ 
\vep(N)  \phi(N) & = &  \frac{1}{2}  ( g_D + i h_D)  
\int_{1}^N  dn \;   \phi(n).  
\nonumber 
\earray

\no
Dividing both eqs,

\beq
\frac{g_D + i h_D}{g_D - i h_D} = 
\frac{ \vep(N)  \phi(N)}{\vep(1)  \phi(1)}, 
\label{b7} 
\eeq

\no and using (\ref{b4}), one finds

\beq
\frac{g_D + i h_D}{g_D - i h_D} = 
{\rm exp} \left( i h_D \int_1^N \frac{d n}{\vep(n)} \right),  
\label{b8} 
\eeq

\no
or equivalently

\beq
\frac{g_I + i h_I}{g_I - i h_I} = 
{\rm exp} \left(\frac{i E}{\hbar}  \int_1^N \frac{d n}{\vep(n)} \right). 
\label{b9} 
\eeq

\no
The eigenenergies are obtained 
taking the log

\beq
N_{\rm I}(E) = 
\frac{E}{2 \pi \hbar}  \int_1^N \frac{d n}{\vep(n)} - 
\frac{\alpha}{\pi}, 
\label{b10} 
\eeq

\no
with

\beq
\alpha = \arc \left( \frac{h_I}{g_I} \right) 
\label{b11} 
\eeq

\no and where $N_I=0, \pm 1, \pm 2, \dots$ label the eigenstates.
In the case $\vep(n) = n, g_I =0$, (\ref{b10}) reduces to

\beq
N_{I}(E) = 
\frac{ E}{2 \pi \hbar}  \log N -  
\frac{1}{2}
\label{b12} 
\eeq

\no
which yields an equally space spectrum symmetric
around zero energy,

\beq
\frac{E_n}{\hbar} = \frac{ 2 \pi }{\log N} ( n + 1/2), 
\qquad n = 0, \pm 1, \dots  
\label{b13} 
\eeq

\no
In the limit $N \rightarrow \infty$, this spectrum
becomes a continuum with constant energy level
density. 
The conclusion is that the regularized 
BK model does not have a discrete spectrum 
whose eigenenergies could be identified
with the Riemann zeros or an approximation to them. 
On the other hand,  eq.(\ref{b12}) 
agrees with the leading 
term in Connes's  formula (\ref{a10}), 
which also diverges with the cutoff $\Lambda$. 

The wave function (\ref{b5}) in the case
$\vep(n) = n$ is

\beq
\psi_E(n) = \frac{A}{n^{1/2 - i E/\hbar}},
\label{b14}
\eeq

\no
which agrees with the BK wave function (\ref{a13}). 
This is the wave function
of a free particle moving in a box of length
$L_N= \log N$. Indeed, let us define
the coordinate

\beq
q \equiv \log n, \qquad 0 \leq q \leq L_N = \log N,  
\label{b15}
\eeq

\no 
and the momentum variable

\beq
k \equiv \frac{E}{\hbar},  \qquad k_n = \frac{2 \pi}{L_N} 
( n + \frac{1}{2} ).  
\label{b16}
\eeq

\no
The scalar product of two eigenstates (\ref{b14})
with quantum numbers $n_1$ and $n_2$ is

\beq
\langle \psi_{E_{n_1}} | \psi_{E_{n_2}} \rangle
= A^2 \int_1^{N} \frac{d n}{n} n^{i (E_{n_1} - E_{n_2})/\hbar }
= A^2 \int_0^{L_N} d q  \;  e^{2 \pi i q (n_1 - n_2)/L_N} =  
\delta_{n_1, n_2},    
\label{b17}
\eeq

\no
where $A= L_N^{-1/2}$ is a normalization constant. 
Notice that the plane waves satisfy
antiperiodic boundary conditions, which can be changed by
choosing $g_I \neq 0$. In particular,  $g_I = \infty$
yields periodic BC's. 

The variable $q$ can be defined for generic 
choices of levels $\vep(n)$

\beq
q(n) = \int_1^n \frac{d n'}{\vep(n')}, \qquad 0 \leq q \leq L_N =  
 \int_1^N \frac{d n'}{\vep(n')}.  
\label{b18}
\eeq

\no
The momenta $k$ is still defined  by (\ref{b16}).
The eigenenergies $E$, or rather momenta $k= E/\hbar$, satisfy 
(\ref{b9}), which can be written as

\beq
e^{ 2 i \alpha} = e^{ i k L_N}.  
\label{b19}
\eeq

\no
This is the quantization condition of a free
particle moving in a box of length $L_N$, with twisted
boundary conditions fixed by $\alpha$. 
We give below a  RG interpretation 
of $q$.

\section{Renormalization Group Analysis}

The Hamiltonian $H_I$ is renormalizable in the sense that
upon integration of the high energy degrees
of freedom, the effective Hamiltonian governing the dynamics
of the remaining variables 
coincides with the original one parameterized by 
 new coupling constants. For a discrete
QM model the RG procedure consists in the Gauss
elimination of the highest energy degree
of freedom and its replacement into
the Schr\"odinger eq. for the other  
variables (one can also eliminate the lowest energy component in which
case the flow goes towards the ultraviolet) \cite{GW}.
Let us describe this process for $H_I$.
We start from eq.(\ref{a25}) written as

\beq
( E_I -  \frac{1}{2} g_I  f_n^2) \psi_n =  
\frac{g_I + i h_I}{2} 
\sum_{m=1}^{n-1}  f_n f_m \psi_m 
+ \frac{g_I - i h_I}{2} 
\sum_{m=n+1}^{N}  f_n f_m \psi_m. 
\label{c1}
\eeq

\no 
Next we eliminate $\psi_N$ 
in terms of $\psi_{n < N}$

\beq
\psi_N = \frac{ g_I + i h_I}{ 2 E_I - g_I f_N^2}
\; \sum_{m=1}^{N-1}  f_N f_m \psi_m. 
\label{c2}
\eeq

\no
Plugging this equation back into (\ref{c1}), 
for $n < N$, one obtains a system of equations
for the variables $\psi_{n < N}$, which is  
identical to the original system (\ref{c1}) except that 
the couplings are 

\barray
g'_I & = &  g_I + \frac{ g_I^2 + h_I^2}{2 E_I 
f_N^{-2} - g_I}, 
\label{c3} \\
h'_I & = &  h_I, 
\label{c4}
\\
f'_n &= & f_n \qquad (n=1, \dots, N-1).
\label{c5}  
\earray

\no
Hence $h_I$  and $f_n$ are RG invariant
couplings, while $g_I$ changes under the RG. 
This is  the reason to add the $g_I$ coupling to
the Hamiltonian $H_I$ since it is generated by the RG. 

The RD model (\ref{a26}) leads to an equation
similar to eq.(\ref{c3}) for $g_D$, where the term
$E_I f_N^{-2} = E_I \vep_N$ is replaced by 
$\vep_N - E_{\rm RD}$, with $E_{\rm RD}$  the energy of
the eigenstate.
To study the low energy degrees of freedom, 
i.e. $|E_{\rm RD}| << \vep_N$,  one can replace 
 $\vep_N - E_{\rm RD}$ by  $\vep_N$, obtaining 
a RG equation for $g_{D}$ which does
not depend on $E_{\rm RD}$. However, 
to analyze eq.(\ref{c3}) in the continuum, 
we cannot eliminate the term $E_I$,
otherwise the scale dependence would dissapear.
We shall use instead the parametrization (\ref{a29})
to express $g_I, h_I, f_N$ in terms of the
couplings $g_D, h_D, \vep_N$ of the associated RD model.
Eq.(\ref{c3}) implies for $g_D$

\beq
g'_D  =   g_D + \frac{ g_D^2 + h_D^2}{2 \vep_N  
- g_D}, 
\label{c6} 
\eeq

\begin{figure}[t!]
\begin{center}
\includegraphics[height= 5 cm,angle= 0]{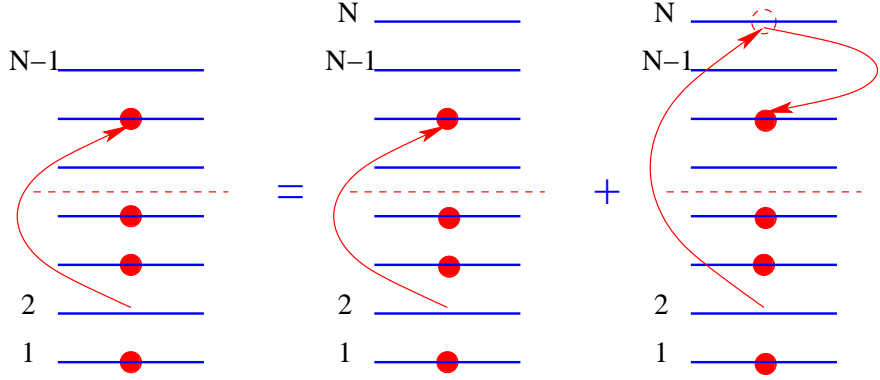}
\end{center}
\caption{
Graphic representation of  eq.(\ref{c6}). The term
$g_D^2 + h_D^2$ arises from the tunneling
of a pair of electrons from a level $n_i <N$ to the highest
level $N$, followed by the decay to a level $n_f < N$. 
The denominator $2 \vep_N - g_D$ is a kinematical
factor.    
}
\label{coupling}
\end{figure}

\no
where $g_D$ and $h_D$ must be regarded
as functions of $E_I$, or $E= 1/E_I$. 
Eq.(\ref{c6}) is the exact RG equation associated to a state
with zero energy in the RD model. This result
is a consequence of the relation
between the $\I$ and RD models explained in section III.
A pictorial representation of eq.(\ref{c6})
is given in fig. \ref{coupling}.

Let's analyze eq.(\ref{c6}) in the continuum limit 
as in the previous section.  For weak
couplings we get

\beq
g_D(N-1) \simeq    g_D(N) + \frac{ g_D^2(N) + h_D^2}{2 \vep_N},  
\label{c7} 
\eeq

\no
where $g_D= g_D(N)$ and  $g'_D= g_D(N-1)$ are
 regarded as a continuous function
of $N$. A gradient expansion of $g_D(N)$ yields 
the differential eq.

\beq
\frac{dg_D}{dN}  = -  \frac{ g_D^2 + h_D^2}{2 \vep_N},  
\label{c8} 
\eeq

\no
which for the uniform model becomes

\beq
\frac{dg_D}{ds}  =   \frac{1}{2}( g_D^2 + h_D^2), 
\label{c9} 
\eeq

\no
where  $N(s)= e^{-s} N$ is the system size at the scale $s$ 
and $N(0)= N$ is the initial size of the system.  
In  general, for a monotonically increasing
function $\vep(n)$, one can define
the scaling variable $s$ as 

\beq
s(n) = \int_{n}^N \;  \frac{d n'}{ \vep(n')},  
\label{c10} 
\eeq

\no
so that eq.(\ref{c8}) takes also the form of eq. 
(\ref{c9}). Comparison of eqs.(\ref{c10})
and (\ref{b18}) yields

\beq
s(n) = L_N - q(n),  
\label{c11} 
\eeq

\no
which relates $s$ to the variable $q$
used in the previous section.  The solution
of (\ref{c9}) is

\beq
g_D(s) = h_D \tan \left( \frac{1}{2} h_D \; s + \arc \left(
\frac{g_D}{h_D} \right) \right),  
\label{c12} 
\eeq

\no
therefore

\beq
g_I(s) = \frac{1}{\hbar}  \tan \left( \frac{E}{2 \hbar} \; s + 
\alpha \right).  
\label{c13} 
\eeq

\no
The coupling $g_I(s)$ is periodic
under the RG with a period

\beq
\lambda_E =  \frac{ 2 \pi}{h_D} = \frac{ 2 \pi \hbar}{E}.  
\label{c14} 
\eeq

\no
In the RD model the RG period
of all the flows is the constant $ 2 \pi/h_D$. However,
(\ref{c14}), yields  a RG period
which depends on the energy
of the state. 
This fact is related to the energy spectrum as
we explain below.

{\bf Russian doll scaling: gapped versus gapless}

Let us consider the $n^{\rm th}$ state $E_n(N)$ 
given by eq.(\ref{b10}),

\beq
n = 
\frac{E_n(N)}{2 \pi \hbar}  \int_1^N \frac{d m}{\vep(m)} - 
\frac{\alpha}{\pi}. 
\label{c15} 
\eeq

\no
After a RG cycle the size is reduced from $N$ to $N(\lambda_n)$
(recall eqs.(\ref{c14}) and (\ref{c10})),

\beq
\lambda_n =  \frac{ 2 \pi \hbar}{E_n(N)} = \int_{N(\lambda_n)}^N
\; \frac{dm}{\vep(m)},  
\label{c16} 
\eeq

\no 
so

\beq
1 = \frac{E_n(N)}{2 \pi \hbar}  \int_{N(\lambda_n)}^N
\; \frac{dm}{\vep(m)}.  
\label{c17} 
\eeq

\no
Splitting the interval of integration $(1,N)$ in 
eq.(\ref{c15}) into the intervals $(1,N(\lambda_n)) \cup 
(N(\lambda_n), N)$, and using (\ref{c17}) one gets

\beq
n-1 = \frac{E_n(N)}{2 \pi \hbar}  \int_1^{N(\lambda_n)}
\; \frac{dm}{\vep(m)} -  \frac{\alpha}{\pi}. 
\label{c18} 
\eeq

\no
This equation has the same form as eq.(\ref{c15})
with the replacements $n \rightarrow n-1$
and  $N  \rightarrow  N(\lambda_n)$, which
implies,

\beq
E_n(N) = E_{n-1}( N(\lambda_n)). 
\label{c19} 
\eeq

\no
In the uniform case $\vep_n = n, g_I = 0$ one has

\beq
\lambda_n = \frac{\log N}{n + 1/2},
\qquad N(\lambda_n) = e^{- \lambda_n} N = N^{\frac{n - 1/2}{n + 1/2}},
\label{c20} 
\eeq

\no
so

\beq
E_n(N) = E_{n-1}(N^{\frac{n - 1/2}{n + 1/2}}), 
\label{c21} 
\eeq

\no 
which can be readily verified using eq.(\ref{b13}). 
In the usual RD model the RG period is the constant $\lambda_D = 2 \pi/h_D$, 
which does not dependent on the energy of the states. The analogue of  scaling
relation (\ref{c21}) is

\beq
E_n(N) = E_{n-1}( e^{- \lambda_D} N). 
\label{c22} 
\eeq

\no
This relation leads to an exponential decaying behaviour of 
the bound state energies

\beq
E_n(N) \sim  N e^{-n \lambda_D},  
\label{c23} 
\eeq

\begin{figure}[t!]
\begin{center}
\includegraphics[height= 6 cm,angle= 0]{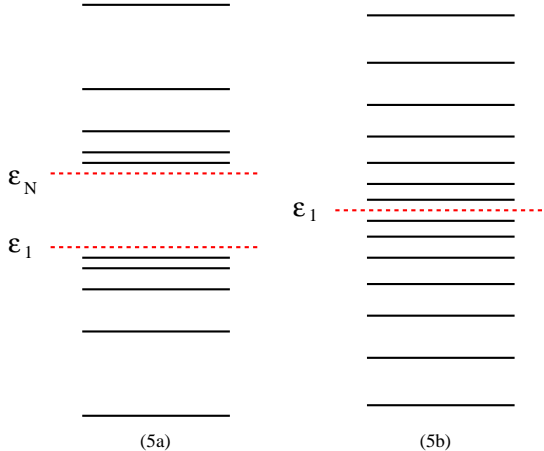}
\end{center}
\caption{
5a) Energy levels of the RD model. The bound states 
lie below the lowest energy $\vep_1$, while the antibound
states lie above the highest energy $\vep_N$. There are also
states within the energy band $(\vep_1, \vep_N)$ (not depicted). 
The scaling is exponential 
(see eq.(\ref{c23})). 5b) Energy levels $E_{I,n} = 1/E_n$ 
of the $\I$ model (see eq.(\ref{b13})).
The energies converge  to the energy  $\vep_1 =0$ 
as $\sim 1/|n|$, from above and below. 
}
\label{levels}
\end{figure}

\no
which is in sharp contrast with the linear decaying of the 
energies $E_{I,n} = 1/E_n$ 
in the $\I$-model (see eq. (\ref{b13})), which reflects 
the gapless nature of its spectrum
(see fig.\ref{levels}).

{\bf Counting RG cycles}

In a finite system the number of RG cycles is of course
finite. In the RD model this number is given 
approximately by \cite{RD1,RD2}

\beq
n_c \sim \frac{h_{D}}{ 2 \pi} \log N.  
\label{c24} 
\eeq

\no
The reason for (\ref{c24})
 is that after each cycle the system size is
reduced by a factor $e^{- \lambda_D}$ where
$\lambda_D = 2 \pi/h_{D}$. Hence, the number of cycles
$n_c$ needed to reduce the system to one site satisfies 
$e^{- n_c \lambda_D} N \sim 1$, which leads to 
eq.(\ref{c24}). This eq. also gives the number
of bound states of the model,   
due to the one-to-one correspondence
between RG cycles and bound states (recall (\ref{c22})). 

For the uniform $\I$-model we expect the formula
(\ref{c24}) to give the number of RG cycles 
for the appropiate value of $h_D$ as a function
of $E_n$, namely

\beq
n_c(E_n)  \sim \frac{E_n}{ 2 \pi \hbar} \log N 
=  n + 1/2.  
\label{c25} 
\eeq

\no
Iterating  eq.(\ref{c21}) one can compute the size of the system after
$n_c$ RG cycles,

\beq
N \rightarrow N^{\frac{n - 1/2}{n + 1/2}} \rightarrow
 N^{\frac{n - 3/2}{n + 1/2}} \rightarrow
\dots  \rightarrow  N^{\frac{n + 1/2 - n_c}{n + 1/2}} \sim  1. 
\label{c26} 
\eeq

\no 
Hence the integer $n$, labelling the state
$E_n$,  coincides roughly with the
number of RG cycles $n_c$ needed to reduce
the system to one site. 
This result is exact to leading order
in $N$,  but as we shall see below there
are finite size corrections which will play an 
important in the discussion.

\section{Exact solution}

The discrete Schr\"odinger eq.(\ref{a25}) can be solved in a very simple
way which parallels the derivation done in the section III.
This solution coincides with  the exact Bethe ansatz solution
of the RD model found in reference \cite{links} for a single 
Cooper pair. There is another derivation of the
exact solution using the RG method of Glazek and Wilson \cite{RD2}. 
The idea is  to keep the energy $E$ in the Gauss elimination
procedure until the system size is one. That gives an exact
equation for the eigenenergies $E$. The renormalizability
of the model makes this procedure doable and that is the
reason for its solvability. 

As explained in section II,  each state with energy $E_I$
of the Hamiltonian (\ref{a19}) coincides
with the zero energy state of the RD Hamiltonian
(\ref{a23}) under the identifications (\ref{a29}). The equation
of that state is given by (\ref{a26}) with $E_{RD} =0$, i.e.

\beq
\vep_n  \phi_n =  \frac{1}{2} 
\sum_{m=1}^N  ( g_D + i h_D \;  \sign(n-m) ) \phi_m. 
\label{d1} 
\eeq

\no
Substracting the eqs. for $\phi_n$ and $\phi_{n+1}$

\beq
\vep_{n+1} \phi_{n+1} - \vep_n \phi_n =
\frac{i h_D}{2}  ( \phi_n + \phi_{n+1}), 
\label{d2} 
\eeq

\no 
gives a discrete differential equation, which is 
a recursion relation for $\phi_{n+1}$ 
as a function of $\vep_n, \vep_{n+1}$ and $\phi_n$,

\beq
\frac{\phi_{n+1}}{\phi_n} =
\frac{ \vep_{n} + i h_D/2 }{ \vep_{n+1} - i h_D/2},
\qquad n=1, \dots, N-1. 
\label{d3} 
\eeq

\no
A gradient expansion of $\phi_n$ and $\vep_n$ in (\ref{d2}) 
reproduces eq.(\ref{b3}). Iterating (\ref{d3}) yields

\beq
\frac{\phi_{N}}{\phi_1} = \prod_{n=1}^{N-1}
\frac{ \vep_{n} + i h_D/2 }{ \vep_{n+1} - i h_D/2}. 
\label{d4} 
\eeq

\no 
 Eqs.(\ref{d3}) is a set of $N-1$ equations while (\ref{d1}) 
contain $N$ equations. Hence there is one more equation to impose.
It is convenient to choose eqs.(\ref{d1}) at the two boundaries,
$n=1$ and $N$ obtaining,

\barray
(\vep_1 - i h_D/2) \phi_1 & = &  \frac{1}{2} ( g_D - i h_D) 
\sum_{m=1}^N \phi_m  \label{d5} \\
(\vep_N + i h_D/2) \phi_N  & = &  \frac{1}{2} ( g_D + i h_D) 
\sum_{m=1}^N \phi_m. \nonumber 
\earray

\no 
Dividing both eqs.

\beq
\frac{ g_D + i h_D}{  g_D - i h_D}
= \frac{(\vep_N + i h_D/2) \;  \phi_N}{ (\vep_1 - i h_D/2) \;  \phi_1},
\label{d6}
\eeq

\no 
and using (\ref{d4}) gives

\beq
\frac{ g_D + i h_D}{  g_D - i h_D}
=  \prod_{n=1}^{N}
\frac{ \vep_{n} + i h_D/2 }{ \vep_{n} - i h_D/2}, 
\label{d7}
\eeq

\no
which can also be written as

\beq
\frac{ g_I + i h_I}{  g_I - i h_I}
=  \prod_{n=1}^{N}
\frac{ \vep_{n} + i E/2 \hbar }{ \vep_{n} - i E/2 \hbar}, 
\label{d8}
\eeq

\no
where $E$ is the energy of the BK Hamiltonian. 
The solutions of eq.(\ref{d8}) are obtained taking the log

\beq
N_I(E) = \frac{1}{2 \pi i} 
\sum_{n=1}^N \log \left( 
\frac{ \vep_{n} + i E/2 \hbar }{ \vep_{n} - i E/2 \hbar}
\right)  - \frac{\alpha}{\pi} \in \Z, 
\label{d9}
\eeq

\no
where $N_I$ is an integer labelling  the eigenstates.

{\bf 1) Completeness of the spectrum.}

Let's write eq.(\ref{d9}) as

\beq
N_I(E) = \frac{1}{\pi} 
\sum_{n=1}^N \arc \left( 
\frac{E}{2 \hbar \vep_n}
\right)  - \frac{\alpha}{\pi}. 
\label{d10}
\eeq

\no
The condition $\vep_n =  f_n^{-2} > 0, \forall n$
 implies that  $N_I(E)$ is a monotonically increasing
function of $E$ varying in the  interval

\beq
- \frac{N}{2} - \frac{\alpha}{\pi}  \leq
 N_I(E) \leq \frac{N}{2} - \frac{\alpha}{\pi}, 
\label{d11}
\eeq

\no
where $- \infty \leq E \leq \infty$. 
There are exactly $N$ integers in the interval
(\ref{d11}) corresponding to all 
the eigenvalues of 
$H_{BK}^{\rm reg}$. For some choices of $N$ and $\alpha$
one may eventually find a solution with $E= \infty$.
If $N$ is even all the solutions are finite. 

{\bf 2) The spectrum of the operator $P$}

A particular  example  of BK Hamiltonian
is given by the choice $v_n = \vep_n = 1, g_I = 0$, 
which corresponds to the operator $P$ (see eq.(\ref{a21})). 
Its  exact eigenvalues, $p_n$, follow inmediately from 
eq.(\ref{d10}),

\beq
p_n = 2 \tan \left[
\frac{\pi}{N} (n + \frac{1}{2} ) \right],
\qquad  - \frac{N}{2} \leq n <  \frac{N}{2}, 
\label{d11-2}
\eeq

\no
which in the continuum limit coincide with  the
momenta of a free particle in a box with antiperiodic
 BC's, $p_n = \frac{2 \pi}{N} (n + 1/2)$.

{\bf 3) Relation with the continuum approximation.}

The power expansion of  the function $\arc$ in (\ref{d10})
is,

\beq
N_I(E) =  \frac{E}{2 \pi \hbar} 
\sum_{n=1}^N 
\frac{1}{\vep_n}
  - \frac{\alpha}{\pi} + O(E^3).  
\label{d12}
\eeq

\no
In the continuum limit

\beq
\sum_{n=1}^N 
\frac{1}{\vep_n} \rightarrow 
\int_{1}^N 
\frac{dn}{\vep(n)}, 
\label{d13}
\eeq

\no
one recovers eq.(\ref{b10}).
The higher order powers of $E$ in (\ref{d10})
are multiplied by terms of the form $\sum 1/\vep_n^{1 + 2m}$
with $m \geq 1$. In the case where $\vep_n = n$
only the sum  $\sum 1/\vep_n \sim \log N$ diverges
with $N$, which controls the large $N$ limit
of  $N_I(E)$.

{\bf 4) Continuum approximation of the exact solution}

The previous discussion suggests to take the  continuum limit
directly in the exact equation (\ref{d10}), namely

\beq
N_I(E) \eqsim  \frac{1}{\pi} 
\int_{1}^N dn  \;  \arc \left( 
\frac{E}{2 \hbar \vep(n)}
\right)  - \frac{\alpha}{\pi}.  
\label{d14}
\eeq

\no
In the uniform case, $\vep(n) = n, \alpha = \pi/2$, and in the limit
$N >> |E|/2 \hbar >> 1$,  eq.(\ref{d14}) becomes

\beq
N_I(E) \eqsim \frac{E}{2 \pi \hbar}
\log N -  \frac{E}{2 \pi \hbar} \left(
\log  \frac{|E|}{2 \hbar}  - 1 \right) + O(1).  
\label{d15}
\eeq

\no
The finite piece of this equation agrees with  Connes's formula
(\ref{a10}), except for a term linear in $E$, namely $E/( 2 \pi) \log \pi$.
The origin of this term, in the Riemann counting formula, can be traced
back from the expression (\ref{a6}), and it is due to a factor 
$\pi^s$ appearing in the functional relation satisfied by $\zeta(s)$.
The leading term $E/(2 \pi) \log E$, as well as another linear term
in $E$, come from the log of the Gamma function in (\ref{a6}).

The divergent terms in eqs.(\ref{d15}) and (\ref{a10})
have a similar form which depend on the respective cutoffs
$N$ and $\Lambda$. Let us suppose for a while 
that both cutoffs are related by the eq.  $N = \pi \Lambda^2$.
Then the divergent and the finite parts 
in both formulas agree. Notice that the $\pi$ factor 
in the latter relation explains the missing factor 
$E/( 2 \pi) \log \pi$ in eq.(\ref{d15}). 
Unfortunately, the relation  $N = \pi \Lambda^2$
does not seem to follow from the counting 
of states in both models. In the $\I$ model
this is given by $N$, while in the Connes
model, at the semiclassical level, 
it would be given by  $\Lambda^2/\pi \;(\hbar = 1)$ 
corresponding to the phase space
$(0 < x < \Lambda, - \Lambda < p < \Lambda)$. 
This comparison gives $N = \Lambda^2/\pi$, 
rather than  $N = \pi \Lambda^2$. 
In any case, the linear term $\propto E$, in both 
counting formulas, depends on the cutoffs and hence
on the particular regularization choosen. 
Apparently, Connes's regularization 
and ours are different. In the next paragraph
with shall consider a zeta function regularization of the model,
which will shed further light on this issue. 

Another interesting point concerns the spectroscopic
interpretation of our results. As we said above they 
are along the lines of Connes's absortion picture. 
However it must be kept in mind that the eigenstates
counted by the smooth part of the 
Riemann formula  are not really missing in the whole spectra
but shifted to higher energies because the interactions. 
This blueshift makes that in a range of energies, say $(0,E)$, 
there are less states than expected from the analysis
of the continuum limit. In this sense, a more
appropiate spectroscopic interpretation of our results
will be in terms of a blueshift of energy levels,
which are then missing in fixed energy intervals.

{\bf 5) Exact solution in the uniform case.}

Eq.(\ref{d10}) can be given an exact
analytic formula in terms of known functions
in the uniform case. We shall 
add  a  {\em zero point} contribution, $a$,
to the energy levels, i.e. 

\beq
\vep_n = n + a, \qquad n = 1,2, \dots, N, 
\label{d16}
\eeq

\no
which does not modify the large $N$ 
properties of $N_I$ discussed above. It is more
convenient to write (\ref{d10}) 
in the product form (\ref{d8}) ($\hbar = 1$)

\beq
e^{ 2 \pi i N_I(E)} = e^{- 2 i \alpha} 
 \prod_{n=1}^{N}
\frac{ \vep_{n} + i E/2 }{ \vep_{n} - i E/2 }.
\label{d17}
\eeq

\no
Inspired by (\ref{d15}), we shall define
the finite part of $N_I$ as

\beq
n_I(E) = \lim_{N \rightarrow \infty}
\left( \frac{E}{2 \pi} \log N - N_I(E) \right),  
\label{d18}
\eeq

\no
where a minus sign has been introduced 
to take care of the relative minus in eq.(\ref{d15}). 
The expression for $n_I(E)$ follows from (\ref{d17}),

\beq
e^{-  2 \pi i n_I(E)} = \lim_{N \rightarrow \infty} 
e^{ - i E \log N} 
 e^{- 2 i \alpha} 
 \prod_{n=1}^{N}
\frac{ \vep_{n} + i E/2  }{ \vep_{n} - i E/2 }. 
\label{d19}
\eeq

\no
Using the eqs.

\barray 
\frac{1}{\Gamma(z)}  & = &  z e^{\gamma z} 
\prod_{n=1}^\infty \left( 1 + \frac{z}{n} \right)
e^{- z/n},  
\label{d20} \\
\gamma & =  & \lim_{N \rightarrow \infty} 
\left( \sum_{n=1}^N  \frac{1}{n} - \log N \right),  
\label{d21} 
\earray

\no
one finds

\beq
\lim_{N \rightarrow \infty} 
e^{ - i E \log N} 
\prod_{n=1}^{N}
\frac{ n + a  + i E/2  }{ n + a - i E/2 }
= \frac{\Gamma \left( 1 + a - i E/2 \right) }{
 \Gamma \left( 1 + a + i E/2 \right)}, 
\label{d22}
\eeq

\no
which plugged in (\ref{d19}), and taking the log, gives

\beq
n_I(E) = \frac{1}{\pi} \log  \Gamma \left( 1 + a + i \frac{E}{2} \right)
+ \frac{\alpha}{\pi}.  
\label{d23}
\eeq

\no
Comparing this equation with the smooth part
of the Riemann's counting formula (\ref{a6}), we deduce 
that the Gamma term in both eqs. agree if we choose

\beq
a= - \frac{3}{4}.  
\label{d24}
\eeq

\no
The energy levels $\vep_n$ (see eq. (\ref{d16}))  corresponding 
to this choice are those of a harmonic oscillator
with zero point energy $1/4$, instead of $1/2$. 
This vacuum energy arises from Neumann BC's at the origin
which select the
even eigenfunctions under parity. In this sense
the RD model can also be seen as a harmonic oscillator with 
Neumann BC's perturbed by the $g_D$ and $h_D$ interactions 
given in eq.(\ref{a23}).  

The asymptotic expansion of (\ref{d23}) gives
the finite part of (\ref{d15}), and it is related to 
the smooth part of the Riemann formula 
(\ref{a6}), for $\alpha = \pi/2$, as

\beq
N_{\rm sm}(E)  =
n_I(E) - \frac{E}{ 2 \pi} \log \pi  +  \frac{1}{2}.  
\label{d25}
\eeq

\no
Notice again the 
term $- E/ 2 \pi \log \pi$  
which, as we said, is regularization dependent. 
To highlight this point,  we shall apply
a zeta function regularization to the expression
(\ref{d19}). Choosing $N= \infty$ one has formally,

\beq
e^{-2 \pi i \hat{n}_I(E)} = 
 e^{- 2 i \alpha} 
 \frac{ \prod_{n=1}^{\infty} \mu ( \vep_{n} + i E/2)}{ 
\prod_{n=1}^{\infty} \mu ( \vep_{n} - i E/2) 
},
\label{d19-2}
\eeq

\no where $-\hat{n}_I(E)$ stands for the regularized value
of $N_I(E)$ and $\mu$ is a regulator.
The infinite products in (\ref{d19-2}) can
be regularized using the zeta function \cite{den}

\beq
 \prod_{n=0}^{\infty} \mu ( n  + z) =  \mu^{1/2 - z} 
\; \left( 
\frac{1}{\sqrt{2 \pi}} \Gamma(z) \right)^{-1}.   
\label{d19-3}
\eeq

\no
The result is

\beq
e^{-2 \pi i \hat{n}_I(E)} = 
 e^{- 2 i \alpha} \mu^{-i E} \; 
\frac{\Gamma \left( 1 + a - i E/2 
\right)}{\Gamma \left( 1 + a + i E/2 \right)},  
\label{d19-4}
\eeq

\no
and

\beq
\hat{n}_I(E) = \frac{1}{\pi} \log  \Gamma \left( 1 + a + i \frac{E}{2} \right)
+ \frac{\alpha}{\pi} + \frac{E}{2 \pi} \log \mu.  
\label{d19-5}
\eeq

\no
Choosing $\mu = 1/\pi$ one gets the factor  
 $- E/ 2 \pi \log \pi$, and the relation (\ref{d25})
is replaced by

\beq
N_{\rm sm}(E)  =
\hat{n}_I(E) +  \frac{1}{2}.  
\label{d19-6}
\eeq

\begin{figure}[t!]
\begin{center}
\includegraphics[height= 7 cm,angle= 0]{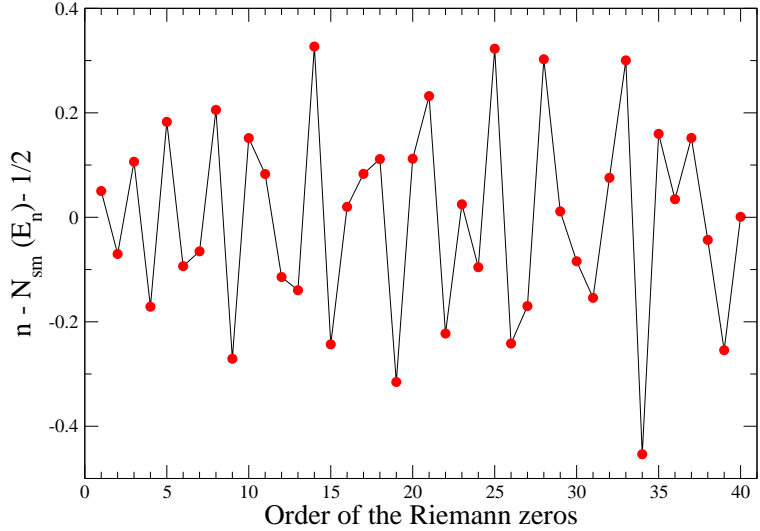}
\end{center}
\caption{
Deviation of the first 40 zeros. The mean and square mean
of $n - N_{\rm sm} - 1/2$ are given by 
 -0.003825  and  0.188102 respectively. 
}
\label{zeros_deviation}
\end{figure}

\no 
There is also a mismatch of $1/2$ due
to the $\alpha$ term, which 
is not important for large values
of $E$ but which can be relevant for small ones. 
We have carried out
a numerical computation to assess the accuracy
of  $N_{\rm sm}(E_n)$ for predicting the position 
of the zeros ($E_n$ is the 
imaginary part of the $n^{\rm th}$ Riemann
zero). The closest  $N_{\rm sm}(E_n)$ comes
to $n$, the better the approximation is. 
The results for the first 10 zeros are 
collected in table 1. They show that  
$N_{\rm sm}(E)+ 1/2 =\hat{n}_I(E) + 1  $
is a much better fit than $N_{\rm sm}(E)$,
suggesting that the zeros are associated to 
complete RG cycles. 
 This result is confirmed 
in figure \ref{zeros_deviation}, where we plot
the difference  $n-N_{\rm sm}(E_n)- 1/2$ 
for the first 40 zeros.  This extra factor
1/2, that improves the location of the zeros,
was also obtained by Berry \cite{B-chaos}
and Badhuri et al. \cite{Bha}
in their respective approaches. 

\begin{center}
\begin{tabular}{|c|c|c|c|c|c|c|c|c|c|c|}
\hline
$n$ & 1 &2  &3  &4 &5 &6 &7 &8 &9 &10   \\
\hline  
$E_n$ & 14.1347 &  21.0220  &  25.0109 & 
  30.4248 & 32.9350 &37.5861  & 40.9187  & 43.3270 & 
48.0051 & 49.7738  \\
$N_{\rm sm}(E_n)$  &  0.4497  &  1.5702 & 2.3936  & 3.6710 &
 4.3172 &5.5935 & 6.5651 & 7.2943 & 8.7708 &  9.3483  \\
 $N_{\rm sm}(E_n)+ \frac{1}{2}$  &  0.9497
 & 2.0702  &  2.8936 & 4.1710 &
  4.8172 & 6.0935 & 7.0651 & 7.7943 & 9.2708 & 9.8483  \\
\hline
\end{tabular}

\vspace{0.5 cm}

Table 1.- Values of the smooth part of the Riemann counting formula,
$N_{\rm sm}(E_n)$ for the first 10 Riemann zeros $E_n$ (eq.(\ref{a6})). 
$N_{\rm sm}(E_n) + 1/2$ is also given for comparison. 
\end{center}

The oscillations  in $n-N_{\rm sm}(E_n)- 1/2$
are randomly distributed due to its relation to 
$\zeta(1/2 + i E)$. This fact implies that
the uniform $\I$-model, eventhough
describes the smooth part of the Riemann's zeros, 
does not explain the origin of their randomness.

To solve this problem we have tried to modify 
the uniform energy levels (\ref{d16}) in several ways. 
The first one is to modify the position
of the energy levels. 
A slight modification of  (\ref{d16})
consist in the addition of a $1/n$ correction,

\beq
\vep_n = n + a + \frac{c}{n}, \qquad n = 1,2, \dots, N,
\label{d26}
\eeq

\no
which does not change the $\log N$ behaviour of $N_I$.
Using (\ref{d19}) one readily finds

\beq
n_I(E) = \frac{1}{\pi} \log  \left[ 
\Gamma \left( 1 + a_+(E)\right)
\;  \Gamma \left( 1 + a_-(E)\right) \right] 
+ \frac{\alpha}{\pi},  
\label{d27}
\eeq

\no
where

\beq
a_\pm(E) = \frac{1}{2} 
\left[ a + i \frac{E}{2} \pm \sqrt{ 
\left(   a + i \frac{E}{2} \right)^2 - 4 c} \right].  
\label{d28}
\eeq

\no
For large values of $E$, eq.(\ref{d27}) agrees with
eq.(\ref{d15}) to order $E \log E$,  but the comparison
with (\ref{a6}) is lost since we have now the product
of two gamma functions instead of one, as in (\ref{d23}). 
Adding higher order corrections to  
(\ref{d26}) of the form $1/n^2, \dots$ 
gives similar results.

Another possibility is to eliminate some 
energy levels in $\vep_n$,  
for example those associated to 
the prime numbers. This choice is suggested by the 
quantum chaos hypothesis according to which
the primitive orbits of the chaotic Hamiltonian
are labelled by the primes \cite{BK1,BK2}. 
The truncation of the prime energy levels, i.e. 
$\vep_p \;$ (with $p$ a prime),  changes the asymptotic
expansion of $N_I$ which behaves as $\log( N/\log N))$
instead of $\log N$. However its finite part, $n_I$, 
does not improve the location of the Riemann zeros. 
These negative results can be understood from the formula
(\ref{d10}) for $N_I$. Indeed, the function $\arc (E/2 \vep_n)$
varies smoothly between 0 and $\pi/2$ in a 
range of $E$ set by  $\vep_n$.  
The same applies for a superposition of those terms with
different energies, making very difficult to obtain 
a random curve interpolating the zeros.  
The uniform energy levels seems to provide  
the best possible approximation within the $\I$ model. 
This fact has lead us to 
a further generalization of this model.

\section{The $\I_\pm$ models}

The couplings $g_{I,D}$ have so far played 
an auxiliary role in the construction of the $\I$ model. They appear
in the renormalization and the exact solution, 
but they did not take part in the dynamics except
for setting the boundary conditions. 
This suggests that $g_{I,D}$ must play a 
more significant role in the dynamics underlying the
Riemann zeros. A natural generalization  is to replace
$g_I$, in the Hamiltonian (\ref{a19}), 
by a generic real symmetric matrix $g^I_{n,m}$, 
however the $h_I$ interaction,
no longer keeps its simple form under the RG,
 becoming  a generic matrix. 
Surprisingly enough, there are two choices of $g^I_{n,m}$, 
which leave the $h_I$ term invariant under the RG. The
corresponding Hamiltonians are given
by,

\beq
\langle n| H_{I_\pm} | m \rangle
= \frac{1}{2}  f_n \; \left( g_{I_\pm}(n,m) + i h_I \;  
\sign(n-m) \right) f_m,  
\label{e1}
\eeq

\no
where the matrix elements $g_{I_\pm}(n,m)$ are defined in terms of a set
of couplings $g_{I,p} \; (p=1, \dots, N)$ as

\beq
g_{I_\pm}(n,m) = g_{I,p}
, \qquad p = 
\left\{ 
\begin{array}{ll}
{\rm max}(n,m) & {\rm for} \;\; \I_+ \\
{\rm min}(n,m) & {\rm for} \;\; \I_- \\
\end{array}. 
\right.
\label{e2}
\eeq

\no

\begin{figure}[t!]
\begin{center}
\includegraphics[height= 6 cm,angle= 0]{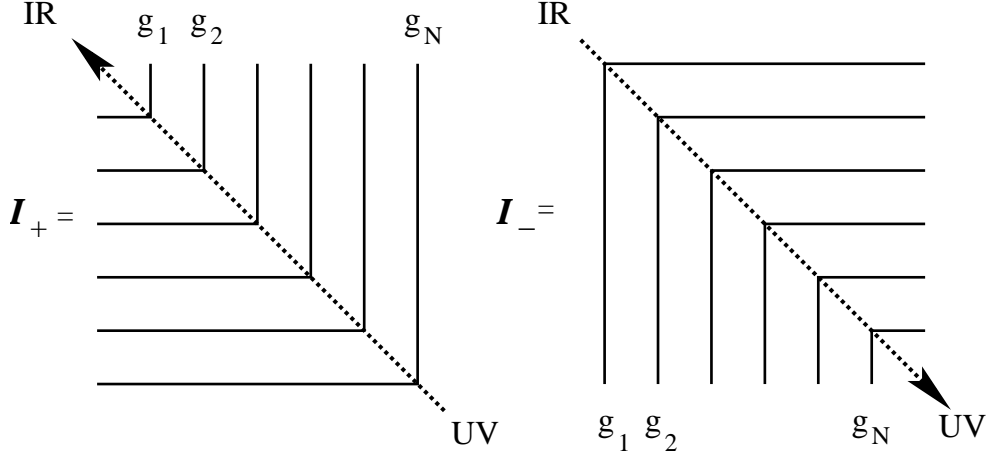}
\end{center}
\caption{
Graphical representation of the coupling matrices $g_{I_\pm}(n,m)$
given in eq.(\ref{e2}). The continuous lines represent common values
of $g_p = g_{I,p}$. The labels $n,m$ (and the energies $\vep_n$)
increase from left to right and from top to bottom. 
The RG in the $\I_+$ model runs towards the infrared (IR) 
and in the  $\I_-$ model towards the ultraviolet (UV). 
}
\label{ipm-model}
\end{figure}

\no 
The structure of  $g_{I_\pm}(n,m)$ 
is displayed in fig. \ref{ipm-model}, which also shows
the direction in which the models are renormalized.  
For $N$ finite the $\I_\pm$  models are related by the 
 transformation $n \rightarrow N-n$ and $h_I \rightarrow
- h_I$, which changes the order of the energy levels $\vep_n$. 
However in the limit $N \rightarrow \infty$ they 
are inequivalent. We shall assume, as usual,
that $\vep_n$ increases with $n$. 

The correspondence between eigenstates
of the $\I$ and the RD  models is also
maintained for $\I_\pm$, defining the associated 
RD$_\pm$ models as

\beq
\langle n| H_{RD_\pm} | m \rangle
= \vep_{n} \;  \delta_{n,m}
- \frac{1}{2}  \; \left( g_{D_\pm}(n,m)  + i h_D \;  \sign(n-m) \right), 
\label{e3}
\eeq

\no where  $g_{D_\pm}(n,m)$ has the form (\ref{e2}), 
with $g_{I,p}$ replaced by $g_{D,p}$.
The relation is established between an eigenstate,
with energy $E_I$, of the $\I_\pm$ models
and the zero energy state, $E_{RD_\pm} = 0$
of the RD$_\pm$, with couplings constants
given by eqs.(\ref{a29}), where the equation
$g_D = g_I/E_I$ is replaced by 
$g_{D,p} = g_{I,p}/E_I$.  
If $g_{I,p}= g_I, \; \forall p$ the $I_\pm$
models reduce to the original $\I$ model and  the same
applies to the RD$_\pm$ models which become the RD one.

From a physical
viewpoint the $\I_+$ and RD$_+$ models
are characterized by the fact that the tunneling 
among two energy levels, induced by the $g$ couplings, 
depends only on the state with the highest level $n$,
i.e. the highest energy one for RD$_+$. 
For the  $\I_-$ and RD$_-$ models, it is the lowest level
that matters. Looking at fig.\ref{coupling} one 
realizes that these models
have to be renormalized in different ways for 
$h_{I,D}$ to be invariant.   
In the RD$_+$ model one has to eliminate
the highest energy mode, $n=N$, which involves
the product $(g_{D,N} + i h_D) (g_{D,N} - i h_D)$
which is real, leaving $h_D$ invariant.
However, in the  RD$_-$ model, the level to be eliminated
is $n=1$, which also  lives $h_D$ fixed. Hence the RG  
flow goes towards the infrared  in the RD$_+$ model, 
and towards the ultraviolet for RD$_-$. 
Of course the RD model can be renormalized in both ways.
We shall not write explicitely the RG eqs. for the
couplings but they can be easily derived using 
the  techniques of section IV.

\subsection{Continuum limit}

Let's  first consider the RD$_+$ model
which, as we said, corresponds to $\I_+$. 
The Sch\"odinger eq.(\ref{e3}), for a state with 
$E_{RD} =0$, is

\beq
\vep_n  \phi_n =  \frac{1}{2} 
\sum_{m=1}^n  g_{D,n} \phi_m + 
\sum_{m=n+1}^N  g_{D,m} \phi_m 
+ i h_D \; \sum_{m=1}^N    \sign(n-m)  \phi_m,  
\label{e4} 
\eeq

\no
and in the continuum 

\beq
\vep(n) \phi(n) = \frac{1}{2} \int_1^n  dm \;   g_D(n)   \;  \phi(m) 
+ \frac{1}{2} \int_n^N  dm \;  g_D(m)  \;  \phi(m) 
+  \frac{i h_D}{2} \int_1^N d m \;  {\rm sign}(n-m) \;  \phi(m).  
\label{e5} 
\eeq

\no
The derivative respect to $n$

\beq
\frac{d}{dn} \left[ 
\vep(n) \phi(n) \right] 
 = i h_D \phi(n) + \frac{1}{2} \frac{d g_D}{dn} 
\int_1^n  dm \;   \phi(m), 
\label{e6} 
\eeq

\no
is an integro-differential equation for $\phi(n)$,  
subject to the boundary condition

\beq
\vep(N) \phi(N) = \frac{1}{2} ( g_D(N) + i h_D)
\int_1^N  dm \;   \phi(m).  
\label{e7} 
\eeq

\no
Eq.(\ref{e6}) can be converted into a second 
order differential equation for the function

\beq
\chi(n) = \int_1^n  dm \;   \phi(m), \qquad \phi(n) 
= \frac{d \chi}{dn},  
\label{e8} 
\eeq

\no
namely

\beq
\frac{d}{dn} \left( 
\vep(n) \frac{ d \chi}{d n} \right) - i h_D 
\frac{ d \chi}{d n} -
 \frac{1}{2} \frac{d g_D}{dn} \chi  = 0.  
\label{e9} 
\eeq

\no
together with the BC's

\beq
\chi(1) = 0, \qquad \vep  \;  \frac{d \chi}{dn} (N) = 
 \frac{1}{2}  ( g_D (N) + i h_D)
\chi(N). 
\label{e10} 
\eeq

\no
Making the change of variables
$n \rightarrow q$ (recall eq.(\ref{b18})), the 
eqs.(\ref{e9}) and (\ref{e10})
become

\beq
\frac{d^2 \chi}{dq^2}  - i h_D 
\frac{ d \chi}{d q} -
 \frac{1}{2} \frac{d g_D}{dq} \chi  = 0, 
\label{e11} 
\eeq

\beq
\chi= 0 \; \; ({\rm at} \;  q=0) , \qquad \frac{d \chi}{dq} = 
 \frac{1}{2}  ( g_D + i h_D)\chi   \; \;    ( {\rm at} \; q= L_N).   
\label{e12} 
\eeq

\no
The second summand in eq.(\ref{e11}) can
be eliminated by the gauge transformation,

\def\chii{\tilde{\chi}}

\beq
\chii = e^{- i h_D q/2} \; \chi,  
\label{e13} 
\eeq

\no
leading to

\beq
\frac{d^2 \chii}{dq^2}  - 
 \frac{1}{2} \frac{d g_D}{dq} \chii
+ \frac{h_D^2}{4} \chii
  = 0,  
\label{e14} 
\eeq

\beq
\chii= 0 \; \; ({\rm at} \; q=0) , \qquad \frac{d \chii}{dq} = 
 \frac{1}{2} g_D \; 
\chii   \; \;  ({\rm at} \; q= L_N).   
\label{e15} 
\eeq

\no
Notice that (\ref{e14}) is the Schr\"odinger
equation of an {\em effective} Hamiltonian

\beq
H_+ = - \frac{d^2}{dq^2} +  V_+(q),
\qquad V_+(q) = 
 \frac{1}{2} \frac{d g_D}{dq}, 
\label{e16} 
\eeq

\no 
corresponding to a potential given by the derivative
of the function $g_D(q)$. The energy of that
state is given by $E_+ = h_D^2/4$. 
It is interesting to rederive 
the results obtained in  previous sections
where $g_D$ is constant, i.e. 
$V_+=0$. The solution of eq.(\ref{e14}) gives
a superposition of plane waves

\beq
\chii = e^{i h_D q/2} + C  e^{-i h_D q/2},
\label{e17} 
\eeq

\no
where $C=-1$ to satisfy the BC at $q=0$ (eq.(\ref{e15})). 
The BC at $q = L_N$ gives the plane wave quantization
(see eqs.(\ref{b8}) and (\ref{b19}))

\beq
e^{i h_D L_N}  = \frac{ g_D + i h_D}{ g_D - i h_D}. 
\label{e18} 
\eeq

The RD$_-$ model can be studied in a similar manner.
We give for completeness the results. 
The function $\chi(n)$ must be defined as

\beq
\chi(n) = \int_n^N  dm \;   \phi(m), \qquad \phi(n) 
= - \frac{d \chi}{dn}.  
\label{e19} 
\eeq

\no
$\chii$ is still given by eq.(\ref{e13}). 
The {\rm effective} Hamiltonian for $\chii$ is

\beq
H_- = - \frac{d^2}{dq^2} +  V_-(q),
\qquad V_-(q) = 
-  \frac{1}{2} \frac{d g_D}{dq},
\label{e20} 
\eeq

\no
and the eigenenergies are given by 
the same formula $E_- = h^2/4$. The BC's also change

\beq
\frac{d \chii}{dq} = 
- \frac{1}{2} g_1 \; 
\chii   \; \;  ({\rm at} \; q= 0),   
\qquad
\chii= 0 \; \; ({\rm at} \; q=L_N). 
\label{e21} 
\eeq

\no 
As in the $\I$ model these results
can be generalized to the discrete case, obtaining 
an exact equation for the eigenvalues of the $\I_\pm$ 
hamiltonians. The analogue of eq.(\ref{d8})
is  a matrix like Bethe equation, related to 
the fact that the wave function $\chii_n$ 
satisfies a second order discrete differential equation.   
The results will be presented elsewhere.

Eqs. (\ref{e11}) and (\ref{e20})  bring the idea of a potential
 $V_{+}$ or $V_-$,  whose scattering theory would produce
the oscillating part of the Riemann counting formula.
Pavlov and Fadeev showed long time ago 
that the zeta function $\zeta(s)$, on the line $Re(s) = 1$,  
 appears in the scattering 
phase shift of particles moving on surfaces
of constant negative curvature \cite{Fadeev,Lax,G2}. This
is a possible direction of research, which is likely to
be  related to quantum chaos.  
On the other hand, our approach is based
on discrete Hamiltonians, so the solution of the 
previous problem is a possible strategy  to find 
the correct choice of the discrete couplings 
$g_n$ leading to the precise location of the Riemann 
zeros.

On more general grounds, let us recall that 
one of the motivations to consider the Hamiltonian $x p$
was the breaking of time reversal symmetry, which
should be related to the GUE statistics of the zeros.
This Hamiltonian gives an accurate
semiclassical description of them but not of their fluctuations.
The $\I_\pm$ models break time reversal, 
but they also break the reversal
of the RG time direction. Could this additional breaking  be 
related to the GUE statistics of the zeros?

There are another interesting questions regarding
the integrability of the models considered in this paper.
The RD model, including its many body version, is 
exactly solvable \`a la Bethe and integrable 
(i.e. infinite number of conserved quantities)
 \cite{links}. Thanks to the 
map $\I$  $\rightarrow$ RD, we have
been able to solve the $\I$ model and in turn the regularized
BK model. Is the $\I$ model, including its
many body version, integrable? It is not obvious how to
generalize the map  $\I$  $\rightarrow$ RD to the many body
case, due to the existence of several rapidity variables.  
The  QM  $\I_\pm$ models are renormalizable and
 exactly solvable, although the analogue
of the Bethe equation has a more complex structure. 
In the RD model
the coupling $g_D$ parameterizes a boundary operator in the 
transfer matrix of an inhomogenous vertex model. The existence
of several couplings $g_{D,n}$ suggests that the corresponding
transfer matrix, if it exists, involves more than one boundary operator.

\section{Conclusions and Prospects}

We have shown in this paper that the Berry and Keating
Hamiltonian can be quantized in a consistent way
on a lattice. This quantization has been
achieved thanks to the relation of the BK Hamiltonian
to two  QM models which have a cyclic RG, 
specially the RD model of superconductivity. 
The later models are renormalizable and have
an exact solution \`a la Bethe, which permits 
a detailed study of their spectrum. 

The first of these models, i.e. the $\I$ model, 
is the inverse of the BK Hamiltonian
plus an extra coupling which, in the BCS framework,
is the standard pairing interaction. This model
has a continuum spectrum  related to 
RG cycles with an energy dependent period. 
The second model is the RD model, which
has bound states related by a scale factor
which depends on a fixed RG period.
These two models are intimately
related, since an eigenstate with energy $E$ 
of the $\I$ model, 
can be mapped into a bound state  at the threshold
of an RD model whose RG period is given by $E$. 
The map, which  is exact, establish an unexpected 
correspondence between the BK model and 
the RD model of superconductivity where the time
reversal symmetry is broken explicitely. 

Using these QM models we have given an spectral 
interpretation of 
the smooth part of the Riemann counting
formula of the non trivial zeros. In the $\I$ model 
it counts  the missing states in the continuum
spectrum below a given energy. The result
depends on the cutoff,  which in our case 
is given by the number of sites. This result
seems to agree with Connes's absortion 
spectral interpretation in the adelic theory,
but the cutoff is  different and a precise comparison 
with ours is not conclusive.
In the RD model, the smooth part of the Riemann
formula gives the number of  missing bound 
states with respect to the leading term which
follows from the scaling properties of the 
cyclic RG. It is a finite size correction
of the Russian doll scaling. In a certain
sense, it can also be seen as an anomaly 
for the discrete RG transformations, i.e. the RG cycles, 
which leave the Hamiltonians invariant.

We have looked for a choice of parameters
of the $\I$ model 
that would  explain 
the oscillating part of the Riemann formula.
In the RD language this means choosing the
energy levels. However 
a numerical study shows 
that the best choice is given 
by equally space energy levels. 
It thus seems that 
the origin of the random
position of the zeros lies
beyond the $\I$ model. 
To explain this randomness 
we propose  a natural
generalization of the $\I$ and RD models 
where the coupling $g$ is replaced
by a set of discrete couplings $g_n$
which depend on the level $n$. 
The guiding principle is renormalizability, 
which in this context means the  RG invariance of 
the time reversal breaking interaction. 
Quite surprisingly,  there are
two models, $\I_\pm$,  satisfying this condition. 
They differ 
in the way the RG procedure is implemented. 
In the $\I_+$ model, the RG runs from the
UV towards the IR, while in the $\I_-$ model  the order
is reversed.

Finally, we have begun to explore
the properties of the $\I_\pm$  models in the continuum limit, 
finding  two QM models with different 
potentials
related to the gradient of the coupling function $g_n$
and different boundary conditions. We suggest
that an appropiate choice of these potentials,
and in turn of the couplings constants
$g_n$, may explain the local fluctuations
of the Riemann zeros. The  $\I_\pm$ models are likely
to be related to other approaches to the Riemann
zeros, specially to quantum chaos. It would be interesting to
investigate that connection.


\noindent {\bf Acknowledgments}. 
I would like to thank A. LeClair for collaborating in the
first stages of this work and for many clarifying discussions. 
I also thank M. Asorey, J. Garc\'{\i}a-Esteve, 
M.A. Mart\'{\i}n-Delgado, G. Mussardo, 
J. Rodr\'{\i}guez-Laguna and J.M. Rom\'an for conversations.  
This work is supported by the 
CICYT of Spain
under the contract BFM2003-05316-C02-01. I also acknowledge the 
EC Commission for financial support via the FP5 Grant
HPRN-CT-2002-00325 and the 
ESF Science Programme INSTANS 2005-2010.

\end{document}